\documentclass{article}
\usepackage{amsmath}
\usepackage{amssymb}
\usepackage{graphicx}
\title{The main theorem of the Galois theory proven with ideas from the first \textit{M\'{e}moire} of Galois.}
\author{Math Dicker, Hoensbroek, the Netherlands}

\begin{document}
\maketitle

\abstract{A proof of the main theorem of the Galois theory is presented using the main theorem of symmetric polynomials. The idea originated from studying the \textit {M\'{e}moire sur les conditions de r\'{e}solubilit\'{e} des \'{e}quations par radicaux} of Galois. The motto \textit{"Read the masters"} pays off.}\\\

\textbf{Introduction}.\\
The Galois theory today is based on automorphisms of a field extension that leave the basic field invariant. Central in this theory is the well-known correspondence theorem usually proven with an Artin-Dedekind lemma. This article will provide a proof of this theorem which is inspired by ideas you can find in the \textit {M\'{e}moire sur les conditions de r\'{e}solubilit\'{e} des \'{e}quations par radicaux.} For translations see \cite{1,2,4}, for the original manuscript see \cite{3}.  The article also refers to the main theorem of the symmetric polynomials wich Galois also often uses. This article therefore shows that the motto \textit{"Read the masters"}$\footnote[1]{Edwards and Tignol support this motto.}$ works.\\\

\begin{figure}[h!]
	\includegraphics[width=\linewidth]{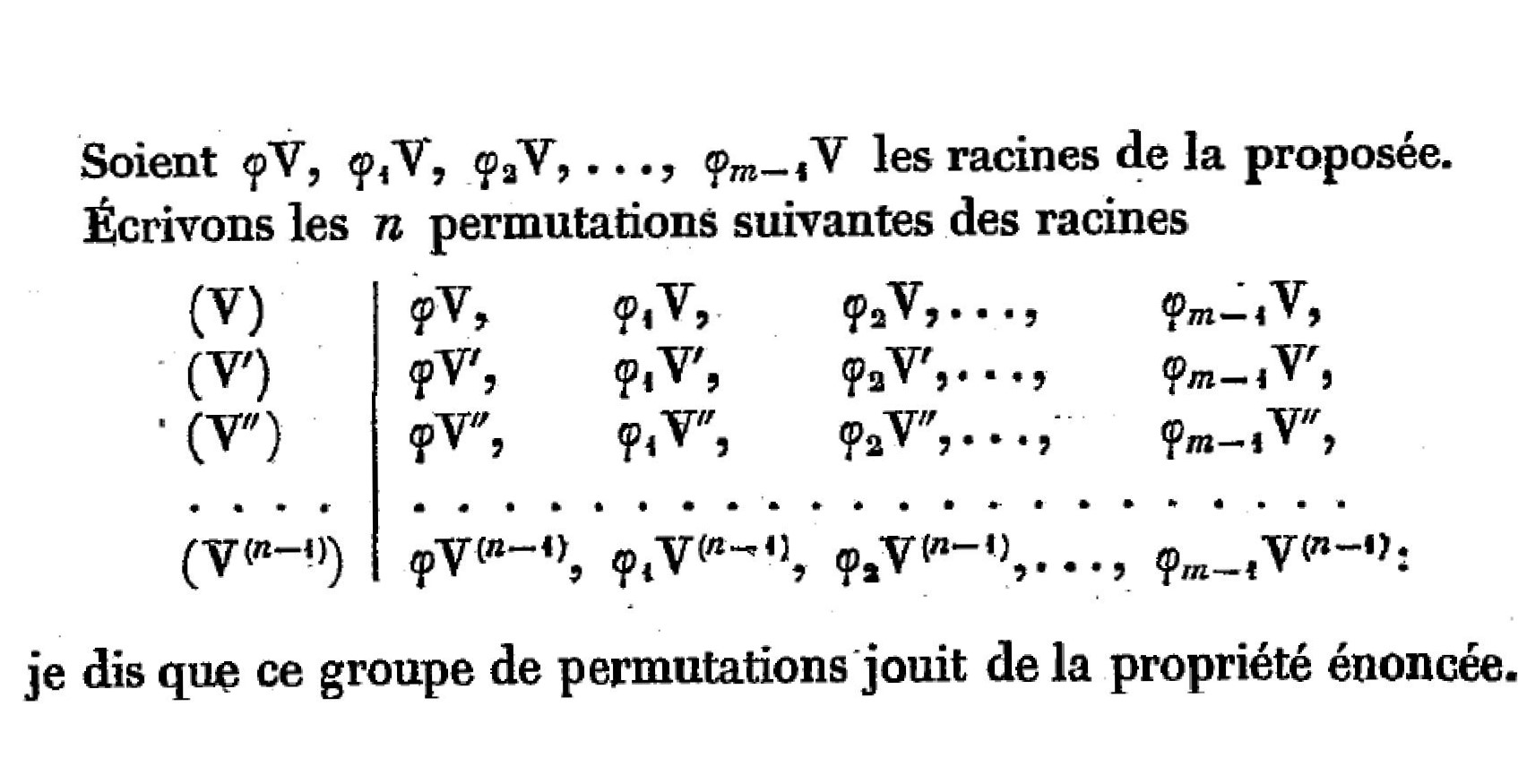}
	\caption{the n$\times$m array published in 1846 in the \textit{Journal de Liouville} in Proposition I.}
\end{figure}

The idea originates as follows: if you look at the n$\times$m array below in proposition I of the \textit {M\'{e}moire} you see that  $\phi$V,$\phi_{1}$V,$\phi_{2}$V,....$\phi_{m-1}$V, \textit{the roots of the given equation}, are listed in different orders in n rows. One entry of the m roots in a row is called an arrangement\footnote[2]{For an excellent explanation  of the concepts arrangement, arrangementgroup, substitution and substitutiongroup see \cite{6}.} of the roots a,b,c,d,.... The n arrangements form a so called arrangementgroup and the corresponding substitutiongroup is the Galois group of F $\subset$ F[a,b,c,d,..] working on a,b,c,d,...  For explanation: an arrangementgroup is a set of arrangements such that if $\alpha$, $\beta$ and $\gamma$$\in$G and $\phi_{\alpha,\beta}$ is the substitution (a bijection on a,b,c,d,....) that transforms $\alpha$ to $\beta$  then $\phi_{\alpha,\beta}$ working on $\gamma$ is also $\in$ G. If G is an arrangementgroup you can make $\Pi$(G)=$\{\phi_{\alpha,\beta} \vert \alpha,\beta \in G \}$, the substitutiongroup of G. This is a group in the modern sense. On the left of the n$\times$m array are the corresponding values from Lemma II; important is that these are all different values. The values V, V$^{'}$, V$^{''}$ ... are all primitive elements for the field F[a,b,c,d,...], this is lemma III. You can make the arrangements in the rows of the n$\times$m array as follows; let every automorfism of the Galois group  transform V from Lemma II. If H is a subgroup of the Galois group we will add a field to H using an arrangementgroup. The idea is this:  let H  be a subgroup of order k of the Galoisgroup with H=$\{h_{1},h_{2},..,h_{k}\}$, choose a random arrangement from the n rows and let $\alpha$ be the corresponding primitive element on the left. We add to H the field  $\{P(h_{1}(\alpha),h_{2}(\alpha),...,h_{k}(\alpha))  \vert$ P a symmetric polynomial $\in  F[x_{1},x_{2},...,x_{k}]\}$. This field appears to be  independent of the choice of the arrangement. If you let H operate on the chosen arrangement this will result in an  arrangementgroup with the same corresponding  primitive elements. If you have an arrangementgroup you can make the substitutiongroup; in this case you get H. But this mapping is not injective; if you read the example before proposition VI where Galois shows that the roots of the quartic are radical, the arrangementgroup of 12 arrangements is partitioned in 3 arrangementgroups of 4 elements. The substitutiongroups of these 3 arrangementgroups are all the four-group of Klein. In case of a fixed arrangement $\epsilon$ then you can prove that there is a one-to-one correspondence between the arrangementgroups containing $\epsilon$ and the subfields of the splitting field mentioned above. We will start from a subgroup of automorphisms in the Galois group and prove the main theorem of the Galois theory in a modern way. The story above is only to illustrate  that is very useful to study the original writings of Galois. I am also convinced that Galois knows everything about symmetric polynomials\footnote[3]{Read the very interesting article \cite{5}}; Galois gives in proposition I a characterization of F that refers strongly to the main theorem of Galois theory. We will  prove this proposition at the end. \\\\

\textbf{The proof}.\\\
A Galois extension can be seen as a splitting field over a basic field F of a polynomial $\in$ F [x] with different roots. We assume that the integers form part of that basic field F in connection with lemma II of Galois. Let R=$\{r_{1},r_{2},...,r_{n}\}$ be the different roots of the aforementioned polynomial then  F $\subset$ F[$r_{1},r_{2},...,r_{n}$] is a Galois extension. Lemma's II and III of Galois  guarantee that there are primitive elements. Choose a primitive element $\omega$ for the  Galois extension of F $\subset$ F[R]. We define a mapping which to every subroup in the Galois group of F $\subset$ F[R]  assigns a field in F[R]; we prove that this mapping is surjective and injective. We also prove that this field is exactly the field referred to in the main theorem proven by E. Artin. Proposition I of in the \textit {M\'{e}moire} is thereby proven.  We  prove the main theorem of the Galois theory using the main theorem of the symmetric polynomials as you will see.\\\

\textbf{The definition af the field L$_{H,\omega}$.}\\\
[0] Let H be a group and H $\subset$ Galois group of F $\subset$ F[R]; assume H = $\{\sigma_{1},\sigma_{2},...,\sigma_{m}\}$ and $\sigma_{i}$ an automorphism of F[R]. Assign to H the field L$_{H,\omega}$= $\{P(\sigma_{1}(\omega),\sigma_{2}(\omega),...,\sigma_{m}(\omega))  \vert$  with P a symmetric polynomial $\in$ $F[x_{1},x_{2},...,x_{m}]\}$; F $\subset$ L$_{H,\omega}$ because  if x $\in$ F use  P$(x_1,x_2,..,x_m)$=x. If x $\notin$ F and x $in$  L$_{H,\omega}$ use the minimum polynomial of x over F to prove that $\frac{1}{x}$ $\in$L$_{H,\omega}$.\\

\textbf{Independence of the choice of the primitive element.}\\\
[1]: If $\theta$ is an other primitive element then L$_{H,\theta}$=L$_{H,\omega}$;   $\theta$=g($\omega$) for some g $\in$ F[x]; Substituting $\theta$=g($\omega$) implies L$_{H,\theta}$$\subset$L$_{H,\omega}$. The other inclusion mutatis mutandis.\\\ 

\textbf{Surjectivity.}
[2]: If L is a subfield of F[R] and g is the minimal polynomial of $\omega$ over L we consider the field  L$_{G,\omega}$ with G the Galois group of L $\subset$ F[R]. Let G=$\{\tau_{1},\tau_{2},...,\tau_{n}\}$ with n equal to the degree of the extension L $\subset$ F[R]. The polynomial g(x)=$\prod_{i=1}^{n} (x-\tau_{i}(\omega))$ is in L[x] and in L$_{G,\omega}$[x]. Using the main theorem of symmetric polynomials and the coefficients of g$\in$L[x]  we know that L$_{G,\omega}$$\subset$L. The minimal polynomial of $\omega$ over L$_{G,\omega}$ divides g and so  [F[R]:L$_{G,\omega}$] $\le$ [F[R]:L]. Therefore also L$_{G,\omega}$=L and the surjectivity is proven.\\\ 

\textbf{Injectivity.}
[3]: [2] has been proven for every subfield of F[R], we apply [2] on  L=L$_{H,\omega}$ with H as in [0]. We are going to use the equality  L=L$_{H,\omega}$=L$_{G,\omega}$ with G as in [2]. Every automorphism $\in$ H is the identity on L=L$_{H,\omega}$, so H$\subset$G, the Galois group of L=L$_{G,\omega}$ $\subset$ F[R]. Therefore m$\le$n. Let f(x)=$\prod_{i=1}^{m} (x-\sigma_{i}(\omega))$ then f $\in$ L[x]=L$_{H,\omega}$[x]; the minimal polynomial g of $\omega$ over  L=L$_{G,\omega}$ in L[x] has degree n and divides f; consequently n$\le$m and so n=m and H=G. The injectivity is proven; the assumption H$_{1}$$\ne$H$_{2}$ and L$_{H_{1},\omega}$=L$_{H_{2},\omega}$ leads to a contradiction because then H$_{1}$=H$_{2}$=the Galois group of (L$_{H_{1},\omega}$=L$_{H_{2},\omega}$) $\subset$ F[R].\\ 

 \textbf {The main theorem of the Galois theory} follows: the mapping which to a subgroup H of the  Galois group of F $\subset$F[R] assigns the field L$_{H,\omega}$ is is a one-to-one correspondence between the subgroups in the Galois group of F $\subset$F[R] and the subfields in F[R]. [F[R]:L$_{H,\omega}$]=$\vert H \vert$ with H the Galois group of  L$_{H,\omega}$ $\subset$ F[R].\\

Yet to prove L$_{H,\omega}$=L$^{H}$ where L$^{H}$=$\{x\vert \sigma(x)=x$ for all $\sigma \in H\}$. Let x $\in$  L$^{H}$; there exists a polynomial h in F[x] with  x=h($\omega$) and so x=$\sigma_{i}(x)=h(\sigma_{i}(\omega))$. Consider the symmetric polynomial P$(x_1,x_2,..,x_m)$ = ($x_1+x_2...+x_m$)/m; x=P(x,x,..,x)= ($h(\sigma_{1}(\omega))+h(\sigma_{2}(\omega))+...+h(\sigma_{m}(\omega)))$/m and consequently x $\in$  L$_{H,\omega}$. That L$_{H,\omega}$ $\subset$ L$^{H}$ is evident. This proves Proposition I in the \textit {M\'{e}moire} of Galois.\\\\

My email address is louis.dicker@ziggo.nl

\end{document}